\documentclass[a4paper]{article}

\usepackage{amsmath, amssymb, amsthm}
\usepackage{lineno,hyperref}
\usepackage{enumerate}

\usepackage{lineno,hyperref}
\usepackage{comment}

\usepackage{mathdots}
\usepackage{enumerate}

\usepackage[all]{xy}

\usepackage{amsthm}

\newtheorem{thm}{Theorem}[section]

\newtheorem{cor}[thm]{Corollary}

\newtheorem{prop}[thm]{Proposition}
\newtheorem{lem}{Lemma}[section]
\newtheorem{rem}{Remark}[section]
\newtheorem{Def}{Definition}[section]
\newtheorem{ex}{Example}

\newtheorem*{thm1}{Theorem 3.2}

\allowdisplaybreaks

\begin{document}
	
		\title{Classification of Toda-type tt*-structures and \(\mathbb{Z}_{n+1}\)-fixed points}
		\author{Tadashi Udagawa}
	    \date{}
	    \maketitle
		
		\begin{abstract}
		We classify Toda-type tt*-structures in terms of the anti-symmetry condition. A Toda-type tt*-structure is a flat bundle whose flatness condition is the tt*-Toda equation (Guest-Its-Lin). We show that the Toda-type tt*-structure can be described as a fixed point of \(e^{\sqrt{-1}\frac{2\pi}{n+1}}\)-multiplication and this ``intrinsic'' description reduces the possibilities of the anti-symmetry condition to only two cases. We give an application to the relation between tt*-Toda equations and representation theory.
		\end{abstract}
		\vspace{10pt}

    {\flushleft{{\it Keywords:} the tt*-Toda -equation, tt*-structure, \(W\)-algebra}

	
\section{Introduction}
In the early 1990's, S. Cecotti and C. Vafa introduced the topological anti-topological fusion (tt*)-structures to describe properties of \(N=2\) supersymmetric field theories \cite{CV1991}, \cite{CV19932}, \cite{CV1993}. In mathematics \cite{D1993}, B. Dubrovin described a tt*-structure as a flat bundle \((E,\eta,g,\Phi)\) with a Higgs field \(\Phi\) over a complex manifold and its flatness condition is called the tt*-equation. The tt*-Toda equation is one of the few explicitly solvable tt*-equations and it is defined by
\begin{equation}
	(w_j)_{t \overline{t}} = e^{w_j-w_{j-1} } - e^{w_{j+1} - w_j },\ \ \ w_j:\mathbb{C}^* \rightarrow \mathbb{R}, \nonumber
\end{equation}
with the conditions
\begin{enumerate}
	\setlength{\leftskip}{8mm}
	\item [(i)] \(w_j = w_{j+n+1}\) (periodicity),
	
	\item [(ii)] \(w_j + w_{l-j-1}=0\) (anti-symmetry condition),
\end{enumerate}
for some \(l \in \mathbb{Z}\) (M. Guest, C.-S. Lin \cite{GL2012}, \cite{GL2014}). In \cite{GIL20151}, \cite{GIL20152}, \cite{GIL2020}, Guest, A. Its and Lin found all global solutions to the tt*-Toda equation satisfying the conditions (i), (ii) and the radial condition \(w_j(t,\overline{t}) = w_j(|t|)\).\vskip\baselineskip

In this paper, we classify tt*-structures whose tt*-equations are the tt*-Toda equations (Toda-type tt*-structures). First, we give an intrinsic description (in other words, frame independent description) of a Toda-type tt*-structure. Physicists considered the ``\(\mathbb{Z}_N\)-symmetry'' on the tt*-structures and this symmetry induces a Toda-type tt*-structure. In mathematics, the \(\mathbb{Z}_N\)-symmetry can be naturally interpreted as a condition with respect to a specific frame (for details see section \ref{s2-2}). This condition depends on the choice of the frame. However, the notion of tt*-structures is independent of the choice of frame and thus, the \(\mathbb{Z}_N\)-symmetry should be described as an intrinsic property. We consider a \(\mathbb{Z}_{n+1}\)-multiplication on tt*-structures
\begin{equation}
	\omega \cdot (E,\eta,g,\Phi) := (E,\eta,g,\omega \Phi),\ \ \ \ \ \ \omega = e^{\sqrt{-1}\frac{2\pi}{n+1}}, \nonumber
\end{equation}
analogous to the \(\mathbb{C}^*\)-action on the moduli of Higgs bundles introduced by Hitchin (see \cite{H1987}). Hitchin used the \(\mathbb{C}^*\)-action to describe properties of the moduli of Higgs bundles. Our main result is given as follows.
\begin{thm1}\label{prop3.1}
	Let \((E,\eta,g,\Phi)\) be a tt*-structure of rank \(n+1\) over \(\mathbb{C}^*\) with a non-vanishing \(\Phi\). If there exists an isomorphism
	\begin{equation}
		\mathcal{T}:(E,\eta,g,\Phi) \stackrel{\sim}{\longrightarrow} \omega \cdot (E,\eta,g,\Phi) = (E,\eta,g,\omega \Phi), \nonumber
	\end{equation}
	of tt*-structures, then there exists a frame of \(E\) such that the tt*-equation with respect to the frame is given by
	\begin{equation}
		(w_j)_{t \overline{t}} = e^{w_j - w_{j-1}} - e^{w_{j+1}-w_j}, \nonumber
	\end{equation}
	with the condition
	\begin{align}
		&w_j + w_{n-j} = 0 \ \ \ {\rm if}\ \ \mathcal{T}^{n+1} = -Id_E, \nonumber \\
		&w_j + w_{n+1-j} = 0 \ \ \ {\rm if}\ \ \mathcal{T}^{n+1} = Id_E, \nonumber
	\end{align}
	for \(j=0,\cdots,n \), where \(w_{-1} = w_n, w_{n+1} = w_0\), i.e. \((E,\eta,g,\Phi)\) is a Toda-type tt*-structure.\\
	We call a fixed point of \(\mathbb{Z}_{n+1}\)-multiplication a \(\mathbb{Z}_{n+1}\)-fixed point.
\end{thm1}\vspace{2mm}

We prove that a fixed point of the \(\mathbb{Z}_{n+1}\)-multiplication in the set of isomorphism classes of tt*-structures gives the tt*-Toda equation with the anti-symmetry condition. From Theorem \ref{prop3.1}, we obtain an intrinsic description for a Toda-type tt*-structure
\begin{align}
	&\left\{\text{Toda type tt*-structure}\right\}
	\ \simeq\ \left\{\mathbb{Z}_{n+1}\text{-fixed point with a non-vanishing}\ \Phi\right\}. \nonumber
\end{align}
The \(l \in \mathbb{Z}\) in the anti-symmetry condition depends on the \(n+1\)-th power of the isomorphism. As a corollary (Corollary \ref{cor3.1}), we obtain a classification up to isomorphism of Toda-type tt*-structures. We classify Toda-type tt*-structures by the value \(l\) of the corresponding tt*-Toda equation.\vskip\baselineskip

Regarding the classification of Toda-type tt*-structures, we give some observations. First, we give the equivalence of the tt*-Toda equations in terms of the anti-symmetry condition. From Theorem \ref{prop3.1} and Corollary \ref{cor3.1}, we can see that every tt*-Toda equation with \(l \in \mathbb{Z}\) is equivalent to the tt*-Toda equation with \(l=0\) or \(l=1\) (Theorem \ref{thm1.2}). Here, ''equivalent`` means that solutions for \(l \in \mathbb{Z}\) can be written by using solutions for \(l = 0\) or \(1\). As a application, we apply Theorem \ref{thm1.2} to the tt*-Toda equation which involves two unknown functions introduced by Guest, Lin \cite{GL2014}. Then, we can reduce the possibilities for unknown functions to just three cases. We also consider more general situations (Corollary \ref{cor4.2}).\vskip\baselineskip

Second, we investigate the relation between radial solutions to the tt*-Toda equations with the anti-symmetry condition \(l \in \mathbb{Z}\) and representations of a vertex algebra which L. Fredrickson and A. Neitzke called the ``universal principal \(W\)-algebra of \(\mathfrak{sl}_{n+1} \mathbb{C}\)'' (section 5 of \cite{FN2021}). Fredrickson, Neitzke observed a correspondence between systems of differential equation of the tt*-Toda equation and representations of  the vertex algebra. For the tt*-Toda equation with the anti-symmetry condition \(l = 0\), Guest and T. Otofuji gave a relation between the tt*-Toda equation and representation theory \cite{GO2022}. In this paper, we consider a specific Toda-type tt*-structure (Proposition \ref{prop4.2}) and we apply their observations on the tt*-Toda equation with the anti-symmetry condition \(l \in \mathbb{Z}\) (Corollary \ref{cor6.6}). \vskip\baselineskip

This paper is organized as follows. In section \ref{s2}, we review the definitions of the tt*-structure and the isomorphism of tt*-structures following \cite{FLY2021}. We define the notion of the \(\mathbb{Z}_{n+1}\)-symmetry and a \(\mathbb{C}^*\)-action on tt*-structures. In section \ref{s3}, we give an intrinsic description of the tt*-Toda equation (Toda-type tt*-structure) and we classify Toda-type tt*-structures in terms of the anti-symmetry condition (Theorem \ref{prop3.1}, Corollary \ref{cor3.1}). In section \ref{s4}, we give some observations on Toda-type tt*-structures. We show that the cases \(l=0\) and \(l=1\) are essential in the anti-symmetry condition (Theorem \ref{thm1.2}). As an application, we consider a tt*-Toda equation which involves two unknown functions introduced by Guest, Lin \cite{GL2014}. We also apply the observation of Fredrickson and Neitzke \cite{FLY2021} to Toda-type tt*-structures. The effective Virasoro central charge can be expressed by the \(l\)-asymptotic data of a specific Toda-type tt*-structure.

\section{Preliminaries}\label{s2}
\subsection{tt*-structures}\label{s2-1}
The tt*-structure is a special case of harmonic bundle and it was introduced by Cecotti and Vafa in physics and was formulated by Dubrovin in mathematics. A tt*-structure over \(\mathbb{C}^*\) is defined as follows.
\begin{Def}[C. Hertling \cite{H2003}, H. Fan, T. Yang, Z. Lan \cite{FLY2021}]\label{def2.1}
	A tt*-structure \((E,\eta,g,\Phi) \) over \(\mathbb{C}^*\) is a holomorphic vector bundle over \(\mathbb{C}^*\) with a holomorphic structure \(\overline{\partial}_E\), a holomorphic nondegenerate symmetric bilinear form \(\eta\), a Hermitian metric \(g \) and a holomorphic \({\rm End}(E) \)-valued 1-form \(\Phi\) such that
	\begin{enumerate}
		\item [(a)] \(\Phi \) is self-adjoint with respect to \(\eta \), \vspace{1mm}
		
		\item [(b)] a complex conjugate-linear involution \(\kappa \) on \(E \) is given by \(g(a,b) = \eta(\kappa(a),b) \) for \(a,b \in \Gamma(E) \), i.e. \(\kappa^2 = Id_E \) and \(\kappa(\mu a) = \overline{\mu} a \) for \(\mu \in \mathbb{C},\ a \in \Gamma(E) \), \vspace{1mm}
		
		\item [(c)] a flat connection \(\nabla^{\lambda} \) is given by
		\begin{equation}
			\nabla^{\lambda} = D + \lambda^{-1} \Phi + \lambda \Phi^{\dagger_g },\ \ \ \lambda \in S^1, \nonumber
		\end{equation}
		where \(D = \partial_E^g + \overline{\partial}_E \) is the Chern connection and \(\Phi^{\dagger_g} \) is the adjoint operator of \(\Phi \) with respect to \(g \).
	\end{enumerate}
	
	Given a tt*-structure \((E,\eta,g,\Phi)\), the flatness condition
	\begin{equation}
		\left[\partial_E^g,\overline{\partial}_E \right] = -\left[\Phi,\Phi^{\dagger_g } \right] = -\left(\Phi \wedge \Phi^{\dagger_g} + \Phi^{\dagger_g} \wedge \Phi\right), \nonumber
	\end{equation}
	is called the tt*-equation.
\end{Def}
\vskip\baselineskip

\begin{rem}
	A tt*-structure can be defined over a Riemann surface \(\Sigma\) by replacing \(\mathbb{C}^*\) in Definition \ref{def2.1} with \(\Sigma\). More generally, it can be defined over a complex manifold (see details in \cite{FLY2021}).
\end{rem}

In this paper, we use the isomorphism of tt*-structures introduced by Fan, Lan and Yang (\cite{FLY2021}).
\begin{Def}\label{def2.2}
	Let \((E_j,\eta_j,g_j.\Phi_j),\ j=1,2\) be two tt*-structures over \(\mathbb{C}^*\). A bundle map \(T:E_1 \rightarrow E_2\) of two holomorphic bundles is called an isomorphism from \((E_1,\eta_1,g_1,\Phi_1) \) to \((E_2,\eta_2,g_2,\Phi_2) \) if \(T \) satisfies
	\begin{enumerate}
		\item [(1)] \(\eta_1(a,b) = \eta_2(T(a),T(b)) \),
		\item [(2)] \(g_1(a,b) = g_2(T(a),T(b)) \),
		\item [(3)] \(T((\Phi_1)_X(a)) = (\Phi_2)_X(T(a)) \),
	\end{enumerate}
	for all \(a,b \in \Gamma(E), X \in K_{\mathbb{C}^*} \). Here, we say the tt*-equations of \((E_1,\eta_1,g_1,\Phi_1)\) and \((E_2,\eta_2,g_2,\Phi_2)\) are equivalent.
\end{Def}
As examples, we review tt*-structures constructed from the sinh-Gordon equation.
\begin{ex}[The sinh-Gordon equation]\label{ex1}
	Given a solution \(w:\mathbb{C}^* \rightarrow \mathbb{R}\) to the sinh-Gordon equation
	\begin{equation}
		w_{t \overline{t}} = e^{2w} - e^{-2w},\ \ \ \ \ t \in \mathbb{C}^*. \nonumber
	\end{equation}
	Let \(E = \mathbb{C}^* \times \mathbb{C}^2\) be a trivial vector bundle and \(\{e_0,e_1\}\) the standard frame of \(E\). We define \(\eta\), \(\Phi\) and \(g\) by
	\begin{align}
		&\eta(e_i,e_j) = \delta_{i,1-j},\ \ \ \Phi(e_0,e_1) = (e_0,e_1)\left(\begin{array}{cc}
			0 & 1 \\
			1 & 0
		\end{array}\right)dt, \nonumber\\
		&g(e_i,e_j) = e^{(-1)^j w}\delta_{ij}, \nonumber
	\end{align}
	then \((E,\eta,g,\Phi)\) is a tt*-structure over \(\mathbb{C}^*\) and its tt*-equation with respect to \(\{e_0,e_1\}\) is the sinh-Gordon equation.\\ \qed
\end{ex}
\begin{ex}\label{ex2}
	Given a solution \(w\) to the sinh-Gordon equation. Let \(E = \mathbb{C}^* \times \mathbb{C}^2\) be a trivial vector bundle and \(\{e_0,e_1\}\) the standard frame of \(E\). We define \(\eta,\Phi\) and \(g\) by
	\begin{align}
		&\eta(e_i,e_j) = \delta_{ij},\ \ \ \Phi(e_0,e_1) = (e_0,e_1)\left(\begin{array}{cc}
			1 & 0 \\
			0 & -1
		\end{array}\right)dt,\ \ t \in \mathbb{C}^*, \nonumber\\
		&g(e_0,e_0) = g(e_1,e_1) = \cosh{(w)},\ \ \ g(e_0,e_1) = \overline{g(e_1,e_0)} = -\sqrt{-1}\sinh{{(w)}}, \nonumber
	\end{align}
	then \((E,\eta,g,\Phi)\) is a tt*-structure over \(\mathbb{C}^*\) and it is isomorphic to the tt*-structure in Example \ref{ex1} as tt*-structures. Thus, the tt*-equation is equivalent to the sinh-Gordon equation.\\
	\qed
\end{ex}
In general, the tt*-equation is highly nonlinear, and thus, it can be solved only in very special cases such as the sinh-Gordon equation (\cite{MTW1977}) and the tt*-Toda equation (\cite{GIL20151}, \cite{GIL20152}, \cite{GIL2020}).\vskip\baselineskip

\subsection{The \(\mathbb{Z}_{n+1}\)-symmetry}\label{s2-2}
In physics, the Landau-Ginzburg theory gives a tt*-structure from a certain holomorphic function (see \cite{C1991}), and the symmetry of the function induces the symmetry to the physical system and the corresponding tt*-equation. For the case of the tt*-Toda equation, we consider the \(\mathbb{Z}_{n+1}\)-symmetry.
\begin{Def}[The \(\mathbb{Z}_{n+1}\)-symmetry]
	Let \((E,\eta,g,\Phi)\) be a tt*-structure over \(\mathbb{C}^*\). We say \((E,\eta,g,\Phi)\) has the \(\mathbb{Z}_{n+1}\)-symmetry if there exists a frame \(\{e_j\}_{j=0}^{n}\) of \(E\) such that \(g\) is invariant under the multiplications
	\begin{equation}
		\{e_j\}_{j=0}^{n} \rightarrow \{\omega^j e_j\}_{j=0}^{n},\ \ \ \ \ \omega = e^{\sqrt{-1}\frac{2\pi}{n+1}}, \nonumber
	\end{equation}
	or equivalently, \(g\) satisfies the condition
	\begin{equation}
		g(\omega^i e_i,\omega^j e_j) = g(e_i,e_j), \nonumber
	\end{equation}
\end{Def}
The tt*-structures in Example \ref{ex1} and Example \ref{ex2} have the \(\mathbb{Z}_2\)-symmetry. Obviously, the standard frame \(e_0,e_1\) gives the \(\mathbb{Z}_2\)-symmetry to the tt*-structure in Example \ref{ex1}, but for the tt*-structure in Example \ref{ex2} we need to choose another frame.
\begin{ex}[The \(\mathbb{Z}_2\)-symmetry in Example \ref{ex2}]
	Let \((E,\eta,g,\Phi)\) be a tt*-structure in Example \ref{ex2}. We define a frame \(\{\tau_0,\tau_1\}\) of \(E\) by
	\begin{equation}
		\tau_0 = \frac{1}{2}\left(e_0 - \sqrt{-1} e_1\right),\ \ \ 
		\tau_1 = \frac{1}{2}\left(e_0 + \sqrt{-1} e_1\right), \nonumber
	\end{equation}
	then \(\{\tau_0,\tau_1\}\) give the \(\mathbb{Z}_2\)-symmetry.\\
	\qed
\end{ex}
Regarding Example \ref{ex2}, if we assume the \(\mathbb{Z}_2\)-symmetry with respect to the standard frame \(e_0,e_1\) then the tt*-equation is Laplace's equation, not the sinh-Gordon equation. When we assume the \(\mathbb{Z}_{n+1}\)-symmetry, the choice of the frame is important to solve the tt*-equation.

\subsection{\(\mathbb{C}^*\)-action and \(\rho\)-action}\label{s2-3}
In the study of harmonic bundles, a \(\mathbb{C}^*\)-action is a useful tool. In this paper, we consider a \(\mathbb{C}^*\)-action on tt*-structures analogous to the \(\mathbb{C}^*\)-action on Higgs bundles introduced by N. Hitchin (see \cite{H1987}).
\begin{Def}[\(\mathbb{C}^*\)-action]
	We define a \(\mathbb{C}^*\)-action on a tt*-structure \((E,\eta,g,\Phi)\) by
	\begin{equation}
		\mu \cdot (E,\eta,g,\Phi) = (E,\eta,g,\mu \Phi),\ \ \ \ \ \mu \in \mathbb{C}^*. \nonumber
	\end{equation}
	In particular, we call the action of \(\mu \in \mathbb{C}^*\) \(\mu\)-action. From the definition of tt*-structures, \((E,\eta,g,\mu \Phi)\) is also a tt*-structure.
\end{Def}
The \(\mathbb{C}^*\)-action on tt*-structures induces a \(\mathbb{C}^*\)-action on the set of isomorphism classes of tt*-structures. In this paper, we consider a multiplication of \(\omega = e^{\sqrt{-1}\frac{2\pi}{n+1}}\) (\(\mathbb{Z}_{n+1}\)-multiplication) on the set of isomorphism classes of tt*-structures.\vskip\baselineskip

We also consider another \(\mathbb{C}^*\)-action on tt*-structures analogous to the \(\mathbb{C}^*\)-action on Higgs bundles introduced by Fredrickson and Neitzke (\cite{FN2021}).
\begin{Def}[\(\rho\)-action]\label{def2.5}
	Let
	\begin{equation}
		\rho:\mathbb{C} \times \mathbb{C}^* \rightarrow \mathbb{C}^*: \rho(z,t) = \rho_z(t) := e^z t\nonumber
	\end{equation}
	We define a \(\rho_z\)-action\((z \in \mathbb{C})\) on a tt*-structure \((E,\eta,g,\Phi)\) by
	\begin{equation}
		\rho_z \cdot (E,\eta,g,\Phi) := (\rho_z^*E,\rho_z^*\eta,\rho_z^*g,e^{-z}\rho_z^*\Phi), \nonumber
	\end{equation}
	where
	\begin{equation}
		\rho_z^*E = \{(t,a) \in \mathbb{C}^* \times E\ | \  \rho_z(t) = \pi(a)\},\ \ \ \ \ (\pi:E \rightarrow \mathbb{C}:\text{projection}), \nonumber
	\end{equation}
	is the pullback bundle by \(\rho_z\) and
	\begin{align}
		&(\rho_z^* \eta)_t((t,a),(t.b)) = \eta_{\rho_z(t)}(a,b),\ \ \ (\rho_z^* g)_t((t,a),(t.b)) = g_{\rho_z(t)}(a,b), \nonumber\\
		&(\rho_z^* \Phi)_X(t,a) = (t,\Phi_{(d\rho_z)_t(X)}(a)), \nonumber
	\end{align}
	for \((t,a), (t,b) \in E,\ X \in T_t\mathbb{C}^*\).
\end{Def}
In this paper, we use the \(\mathbb{Z}_{n+1}\)-multiplication and \(\rho_z\)-action (\(z \in \sqrt{-1}\mathbb{R}\)) on tt*-structures.

\section{Classification of Toda-type tt*-structures}\label{s3}
In this section, we classify tt*-structures whose tt*-equations consist of the tt*-Toda equation
\begin{equation}
	(w_j)_{t\overline{t}} = e^{w_j-w_{j-1}} - e^{w_{j+1}-w_j},\ \ \ w_j:\mathbb{C}^* \rightarrow \mathbb{R}, \nonumber
\end{equation}
with the conditions
\begin{enumerate}
	\setlength{\leftskip}{8mm}
	\item [(i)] \(w_j = w_{j+n+1}\) (periodicity),
	
	\item [(ii)] \(w_j + w_{l-j-1}=0\) (anti-symmetry condition),
\end{enumerate}
for some \(l \in \mathbb{Z}\). Here, we do not assume the radial condition.\vskip\baselineskip

In section \ref{s3-1}, we describe the tt*-Toda equation as a tt*-structure (Toda-type tt*-structure) under the assumption of the \(\mathbb{Z}_{n+1}\)-symmetry. However, the \(\mathbb{Z}_{n+1}\)-symmetry depends on the choice of the frame. In section \ref{s3-2}, we consider a fixed point of \(\mathbb{Z}_{n+1}\)-multiplication instead of the \(\mathbb{Z}_{n+1}\)-symmetry. This gives an intrinsic point of view of the tt*-Toda equation with the anti-symmetry condition. In section \ref{s3-3}, we classify Toda-type tt*-structures in terms of the anti-symmetry condition.

\subsection{Toda-type tt*-structures}\label{s3-1}
We construct a tt*-structure whose tt*-equation is the tt*-Toda equation by using the \(\mathbb{Z}_{n+1}\)-symmetry. Let \(E = \mathbb{C}^* \times \mathbb{C}^{n+1}\) be a trivial vector bundle and \(\{e_j \}_{j=0}^{n} \) be the standard frame of \(E \). For \(l \in \{0,\cdots,n\}\), we consider the tt*-equation \((E,\eta,g,\Phi)\) for
\begin{equation}
	\eta(e_i,e_j) = \left\{\begin{array}{ll}
		\delta_{i,l-1-j} & {\rm if}\ 0 \le i \le l-1,\\
		\delta_{i,n+l-j} & {\rm if}\ l \le i \le n,
	\end{array}\right. \label{1} 
\end{equation}
and
\begin{equation}
	\Phi (e_0,\cdots,e_n) = (e_0,\cdots,e_n) \left(\begin{array}{cccc}
		& & & 1 \\
		1 & & &   \\
		& \ddots & &  \\
		& & 1 &
	\end{array}\right)dt,\ \ \ t \in \mathbb{C}^*. \label{2} 
\end{equation}
Here, we assume the \(\mathbb{Z}_{n+1}\)-symmetry with respect to the standard frame. Then the \(\mathbb{Z}_{n+1}\)-symmetry implies that \(g(e_i,e_j) = e^{w_i} \delta_{ij} \) for some functions \(w_i: \mathbb{C}^* \rightarrow \mathbb{R}\).

\begin{prop}
	This \((E,\eta,g,\Phi)\) is a tt*-structure if and only if \(\{w_j\}_{j=0}^{n}\) is a solution to the tt*-Toda equation.
\end{prop}
\begin{proof}
	It is enough to show that the tt*-equation is the tt*-Toda equation. From the condition \(\kappa^2 = Id_E \) the \(\{w_j\}_{j=0}^{n} \) satisfy the anti-symmetry condition
	\begin{equation}
		\left\{\begin{array}{l}
			w_0 + w_{l-1} = 0, w_1 + w_{l-2} = 0, \cdots \\
			w_l + w_n = 0, w_{l+1} + w_{n-1} = 0, \cdots,
		\end{array}\right. \nonumber
	\end{equation}
	where \(w_{n+1+j}=w_j\ (j \in \mathbb{Z})\). From the tt*-equation, the \(\{w_j \}_{j=0}^{n} \) satisfy
	\begin{equation}
		(w_j)_{t \overline{t}} = e^{w_j-w_{j-1} } - e^{w_{j+1} - w_j },\ \ \ j = 0,\cdots,n, \nonumber
	\end{equation}
	where \(w_{-1} = w_n, w_{n+1} = w_0\). This is the tt*-Toda equation.
\end{proof}

The existence of a global solution to the tt*-Toda equation on \(\mathbb{C}^*\) was proved by Guest, Its and Lin from the view point of the p.d.e. theory and isomonodromy theory (see \cite{GIL20151}, \cite{GIL20152}, \cite{GIL2020}) and thus, we obtain a tt*-structure whose tt*-equation is the tt*-Toda equation.\vskip\baselineskip

We obtain the description of the tt*-Toda equation as a tt*-structure. However, we assume the \(\mathbb{Z}_{n+1}\)-symmetry to obtain the tt*-Toda equation. The \(\mathbb{Z}_{n+1}\)-symmetry is expressed by using a frame \(\{e_j\}_{j=0}^{n}\), i.e. it depends on the frame. On the other hand, the notion of tt*-structures is an intrinsic property. In the following section, we describe the tt*-Toda equation as a tt*-structure without using the frame. As a preparation, we define a Toda-type tt*-structure.\vskip\baselineskip

\begin{Def}\label{def3.1}
	Let \((E,\eta,g,\Phi)\) be a tt*-structure over \(\mathbb{C}^*\).\\
	We call \((E,\eta,g,\Phi)\) a Toda-type tt*-structure if there exists a frame \(\{e_j\}_{j=0}^{n}\) of \(E\) such that \(\eta, \Phi\) satisfy (\ref{1}), (\ref{2}) with respect to \(\{e_j\}_{j=0}^{n}\) and \(g(e_i,e_j) = e^{w_j}\delta_{ij}\) for a solution \(\{w_j\}_{j=0}^{n}\) to the tt*-Toda equation.
\end{Def}

We show that a Toda-type tt*-structure is a fixed point of the \(\mathbb{Z}_{n+1}\)-multiplication on the set of isomorphism classes of tt*-structures.
\begin{lem}\label{lem3.1}
	Let \((E,\eta,g,\Phi)\) be a Toda-type tt*-structure, then \((E,\eta,g,\Phi)\) is isomorphic to \((E,\eta,g,\omega \Phi)\) as a tt*-structure.
\end{lem}
\begin{proof}
	Let \(\{e_j\}_{j=0}^{n}\) be the frame in Definition \ref{def3.1}. We define an automorphism \(\mathcal{T}: \Gamma(E) \rightarrow \Gamma(E)\) by \(\mathcal{T}(e_j) = \omega^{j+\frac{1}{2}(1-l)}e_j\). Then we have \(\eta(\mathcal{T}(e_i),\mathcal{T}(e_j)) = \eta(e_i,e_j)\), \(g(\mathcal{T}(e_i),\mathcal{T}(e_j)) = g(e_i,e_j)\) and \(\mathcal{T}(\Phi(e_j)) = \omega \Phi(\mathcal{T}(e_j))\). Thus, we obtain the stated result.
\end{proof}
In the following section, we show the converse, i.e. a fixed point of the \(\mathbb{Z}_{n+1}\)-multiplication gives a Toda-type tt*-structure.

\subsection{An intrinsic description of the tt*-Toda equation}\label{s3-2}
In this section, we consider the \(\mathbb{Z}_{n+1}\)-multiplication on tt*-structures instead of the \(\mathbb{Z}_{n+1}\)-symmetry. In the set of isomorphism classes of tt*-structures,
\begin{align}
	&\{\text{isomorphism class of a Toda-type tt*-structure}\} \nonumber\\
	&\ \ \ \ \ \ \ \ \ \ \subset \left\{\text{fixed point of the}\ \mathbb{Z}_{n+1} \text{-multiplication}\right\}. \nonumber
\end{align}
We show the converse. First, we consider a non-vanishing \(\Phi\).
\begin{lem}\label{lem3.2}
	Let \((E,\eta,g,\Phi)\) be a tt*-structure of rank \(n+1\) over \(\mathbb{C}^*\) with a non-vanishing \(\Phi\) and \(\omega = e^{\sqrt{-1}\frac{2\pi}{n+1}}\). If there exists an isomorphism
	\begin{equation}
		\mathcal{T}:(E,\eta,g,\Phi) \stackrel{\sim}{\longrightarrow} \omega \cdot (E,\eta,g,\Phi) = (E,\eta,g,\omega \Phi), \nonumber
	\end{equation}
	then the eigenvalues of \(\Phi\) are \(1dt,\omega dt,\cdots,\omega^n dt\), where \(t\) is the coordinate of \(\mathbb{C}^*\).
\end{lem}
\begin{proof}
	Let \(z\) be the coordinate of \(\mathbb{C}^*\) and \(u_0 dz,\cdots,u_n dz\) be the eigenvalues of \(\Phi\). Since \(\Phi \neq 0\), there exists a \(j \in \{0,\cdots,n\}\) such that \(u_j \neq 0\). Let \(l_j\) be an eigenvalue frame of \(E\) (not necessary global) and \(\mathcal{T}(l_j) = \sum_{k=0}^{n}T_{kj}l_k\). From \(\mathcal{T}\circ \Phi = (\omega \Phi) \circ \mathcal{T}\), we have \((u_j - \omega u_k)T_{kj} = 0\) for all \(j,k\). Since \(\mathcal{T}\) is isomorphism, for all \(j\) there exists a \(k_j\) such that \(T_{k_j,j} \neq 0\) and \(k_j \neq k_i \ (j\neq i)\). Thus, we have \(\omega u_j = u_{\sigma(j)}\) for some \(\sigma \in \mathfrak{S}_{n+1}\). Furthermore, we can choose the order of \(u_j\)'s so that \(u_j = \omega u_{j-1} = \cdots = \omega^j u_0\) and \(u_0 \neq 0\). By choosing the coordinate \(t = u_0 z\) of \(\mathbb{C}^*\), we obtain the stated result.
\end{proof}
We consider the \(\mathbb{Z}_{n+1}\)-multiplication on tt*-structures with a non-vanishing \(\Phi\). We prove that a fixed point of \(\mathbb{Z}_{n+1}\)-multiplication gives a Toda-type tt*-structure. We show the following lemma.
\begin{lem}\label{lem3.3}
	Let \((E,\eta,g,\Phi)\) be a tt*-structure of rank \(n+1\) over \(\mathbb{C}^*\) with a non-vanishing \(\Phi\). If there exists an isomorphism
	\begin{equation}
		\mathcal{T}:(E,\eta,g,\Phi) \stackrel{\sim}{\longrightarrow} \omega \cdot (E,\eta,g,\Phi) = (E,\eta,g,\omega \Phi), \nonumber
	\end{equation}
	of tt*-structures, then there exists a frame \(\tau = (\tau_0,\cdots,\tau_n)\) of \(E\) on \(\mathbb{C}^*\) such that \begin{equation}
		\mathcal{T}(\tau) = \tau \cdot \left(\begin{array}{cccc}
			& 1 & & \\
			& & \ddots & \\
			& & & 1 \\
			\varepsilon & & &
		\end{array}\right),\ \ \ \varepsilon \in \{-1,1\}. \nonumber
	\end{equation}
	In particular, we obtain \(\mathcal{T}^{n+1} = \varepsilon \cdot Id_E\).
\end{lem}
\begin{proof}
	From Lemma \ref{lem3.2}, the eigenvalues of \(\Phi\) is \(1 dt,\omega dt, \cdots,\omega^n dt\). From Theorem 10.34 of \cite{L2013} (J. Lee, 2013), \(E\) can be split into \(E = L_0 \oplus \cdots \oplus L_n\), where \(L_j\) is a holomorphic subbundle of \(E\) such that \(\Phi(s) = \omega^j dt \cdot s\) for all local section \(s\) of \(L_j\). From Theorem 30.4 of \cite{F1991} (Forster, 1991), \(L_j\) is holomorphically trivial and then there exists a holomorphic frame \(\tilde{\tau} = (\tilde{\tau}_0,\cdots,\tilde{\tau}_n) \) of \(E\) on \(\mathbb{C}^*\) such that \(\Phi(\tilde{\tau}_j) = w^jdt \cdot \tilde{\tau}_j\) and \(\eta(\tilde{\tau}_i,\tilde{\tau}_j) = \delta_{ij}\). Thus, the proof is given by a linear algebra calculation.\vskip\baselineskip
	
	Let \(\mathcal{T}(\tilde{\tau}_j) = \sum_{k=0}^{n}T_{kj}\tilde{\tau}_k\), then from
	\begin{equation}
		\omega \Phi(\mathcal{T}(\tilde{\tau}_j)) = \mathcal{T}(\Phi(\tilde{\tau}_j)) = \omega^j \mathcal{T}(\tilde{\tau}_j)dt, \nonumber
	\end{equation}
	we have \((\omega^{k+1} - w^j)T_{kj} = 0\) for all \(j,k\) and thus
	\begin{equation}
		\mathcal{T}(\tilde{\tau}_0,\cdots,\tilde{\tau}_n) = (\tilde{\tau}_0,\cdots,\tilde{\tau}_n) \left(\begin{array}{cccc}
			& \varepsilon_1 & & \\
			& & \ddots & \\
			& & & \varepsilon_n \\
			\varepsilon_0 & & &
		\end{array}\right),\ \ \ \epsilon_j \neq 0. \nonumber
	\end{equation}
	It follows from \(\eta(\mathcal{T}(\tau_j),\mathcal{T}(\tau_l)) = \eta(\tau_j,\tau_l)\) that \(\varepsilon_j^2 = 1\). Put \(\tau = (\tau_0,\cdots,\tau_n) = (\tilde{\tau}_0,\cdots,\tilde{\tau}_n) \cdot {\rm diag}(\varepsilon_0,\varepsilon_0 \varepsilon_1,\cdots,\varepsilon_0 \cdots \varepsilon_n)\), then we obtain
	\begin{equation}
		\mathcal{T}(\tau_0,\cdots,\tau_n) = (\tau_0,\cdots,\tau_n)\left(\begin{array}{cccc}
			& 1 & & \\
			& & \ddots & \\
			& & & 1 \\
			\varepsilon & & & 
		\end{array}\right) := \tau \cdot T, \nonumber
	\end{equation}
	where \(\epsilon \in \{ -1,1\}\).
\end{proof}

From Lemma \ref{lem3.3}, we obtain the following theorem.
\begin{thm}\label{thm3.2}
	Let \((E,\eta,g,\Phi)\) be a tt*-structure of rank \(n+1\) over \(\mathbb{C}^*\) with a non-vanishing \(\Phi\). If there exists an isomorphism
	\begin{equation}
		\mathcal{T}:(E,\eta,g,\Phi) \stackrel{\sim}{\longrightarrow} \omega \cdot (E,\eta,g,\Phi) = (E,\eta,g,\omega \Phi), \nonumber
	\end{equation}
	of tt*-structures, then there exists a frame of \(E\) such that the tt*-equation with respect to the frame is given by
	\begin{equation}
		(w_j)_{t \overline{t}} = e^{w_j - w_{j-1}} - e^{w_{j+1}-w_j}, \nonumber
	\end{equation}
	with the condition
	\begin{align}
		&w_j + w_{n-j} = 0 \ \ \ {\rm if}\ \ \mathcal{T}^{n+1} = -Id_E, \nonumber \\
		&w_j + w_{n+1-j} = 0 \ \ \ {\rm if}\ \ \mathcal{T}^{n+1} = Id_E, \nonumber
	\end{align}
	for \(j=0,\cdots,n \), where \(w_{-1} = w_n, w_{n+1} = w_0\), i.e. \((E,\eta,g,\Phi)\) is a Toda-type tt*-structure.\\
	We call a fixed point of \(\mathbb{Z}_{n+1}\)-multiplication a \(\mathbb{Z}_{n+1}\)-fixed point.
\end{thm}
\begin{proof}
	Let \(\tau = (\tau_0,\cdots,\tau_n)\) be a frame of \(E\) on \(\mathbb{C}^*\) in Lemma \ref{lem3.3} and \(G, T\) the representation matrices of \(g, \mathcal{T}\) with respect to \(\tau\) respectively. From the condition \(g(\mathcal{T}(\tau_j),\mathcal{T}(\tau_l)) = g(\tau_j,\tau_l) \), we have \(TG = GT\). Since \(\tau\) is a global frame of \(E\), it is enough to give the proof by a linear algebra calculation.\vskip\baselineskip
	
	When \(\mathcal{T}^{n+1} = -Id_E \), we have \(\varepsilon = -1 \). Let
	\begin{equation}
		L = \frac{1}{\sqrt{n+1}}\left(\begin{array}{cccc}
			1 & & & \\
			& \omega^{-\frac{1}{2} } & & \\
			& & \ddots & \\
			& & & \omega^{-\frac{n}{2}}
		\end{array}\right)\left(\begin{array}{cccc}
			1 & 1 & \cdots & 1 \\
			1 & \omega & \cdots & \omega^n \\
			\vdots & \vdots & \cdots & \vdots \\
			1 & \omega^n& \cdots & \omega^{n^2}
		\end{array}\right), \nonumber
	\end{equation}
	then we have \(L^{-1}TL = \omega^{-\frac{1}{2}}{\rm diag}(1,\omega,\cdots,\omega^n)\) and then
	\begin{equation}
		L^{-1} \left(\begin{array}{cccc}
			1 & & & \\
			& \omega & & \\
			& & \ddots & \\
			& & & \omega^n
		\end{array}\right) L = \left(\begin{array}{cccc}
			& & & 1 \\
			1 & & & \\
			& \ddots & & \\
			& & 1 &
		\end{array}\right). \nonumber
	\end{equation}
	Since \(TG = GT \), we obtain \(A = L^{-1}GL = {\rm diag}(e^{w_1}\cdots,e^{w_n},e^{w_0})\), where \(w_j:U \rightarrow \mathbb{R}\). Put
	\begin{equation}
		e = (e_0,\cdots,e_n) = \tau \cdot L \left(\begin{array}{cccc}
			& & & 1 \\
			1 & & & \\
			& \ddots & & \\
			& & 1 &
		\end{array}\right), \nonumber
	\end{equation}
	then \(e\) is a holomorphic frame of \(E\) such that
	\begin{equation}
		\left(\eta(e_j,e_l) \right) = \left(\begin{array}{ccc}
			& & 1 \\
			& \iddots & \\
			1 & &
		\end{array} \right),\ \ \ \left(g(e_j,e_l)\right) = \left(\begin{array}{ccc}
			e^{w_0} & & \\
			& \ddots & \\
			& & e^{w_n}
		\end{array}\right), \nonumber
	\end{equation}
	\begin{equation}
		\Phi(e_0,\cdots,e_n ) = (e_0,\cdots,e_n)\left(\begin{array}{cccc}
			& & & 1 \\
			1 & & & \\
			& \ddots & & \\
			& & & 1
		\end{array} \right)dt. \nonumber
	\end{equation}
	From \(\kappa^2 = Id_E\), the tt*-equation for \((E,\eta,g,\Phi) \) is equivalent to the tt*-Toda equation
	\begin{equation}
		(w_j)_{t \overline{t}} = e^{w_j - w_{j-1}} - e^{w_{j+1}-w_j}, \nonumber
	\end{equation}
	with the condition
	\begin{equation}
		w_j + w_{n-j} = 0. \nonumber
	\end{equation}\vskip\baselineskip
	
	When \(\mathcal{T}^{n+1} = Id_E \), we have \(\epsilon = 1 \). Let
	\begin{equation}
		L = \frac{1}{\sqrt{n+1}}\left(\begin{array}{cccc}
			1 & 1 & \cdots & 1 \\
			1 & \omega & \cdots & \omega^n\\
			\vdots & \vdots & \cdots & \vdots \\
			1 & \omega^n & \cdots & \omega^{n^2}
		\end{array} \right), \nonumber
	\end{equation}
	and put \(e = (e_0,\cdots,e_n) = \tau \cdot L\), then \(e \) is a holomorphic frame of \(E \) such that
	\begin{equation}
		\left(\eta(e_j,e_l) \right) = \left(\begin{array}{cccc}
			1 & & & \\
			& & & 1 \\
			& &\iddots & \\
			& 1 & &
		\end{array} \right),\ \ \ \left(g(e_j,e_l)\right) = \left(\begin{array}{ccc}
			e^{w_0} & & \\
			& \ddots & \\
			& & e^{w_n}
		\end{array}\right), \nonumber
	\end{equation}
	\begin{equation}
		\Phi(e_0,\cdots,e_n ) = (e_0,\cdots,e_n)\left(\begin{array}{cccc}
			& & & 1 \\
			1 & & & \\
			& \ddots & & \\
			& & & 1
		\end{array} \right)dt. \nonumber
	\end{equation}
	From \(\kappa^2 = Id_E \), the tt*-equation for \((E,\eta,g,\Phi) \) is equivalent to the tt*-Toda equation
	\begin{equation}
		(w_j)_{t \overline{t}} = e^{w_j - w_{j-1}} - e^{w_{j+1}-w_j}, \nonumber
	\end{equation}
	with the condition
	\begin{equation}
		w_j + w_{n+1-j} = 0, \nonumber
	\end{equation}
	where \(w_{n+1} = w_0\).
\end{proof}
The \(\mathbb{Z}_{n+1}\)-fixed point has a tt*-structure with the \(\mathbb{Z}_{n+1}\)-symmetry and its tt*-equation is the tt*-Toda equation with the anti-symmetry condition for \(l=0,1\). The choice of the \(l=0,1\) 
depends on the \(n+1\)-th power of isomorphism \(\mathcal{T}\). This gives an intrinsic description of the Toda-type tt*-structure.

\subsection{Classification of Toda-type tt*-structures}\label{s3-3}
From Theorem \ref{prop3.1}, we obtain an intrinsic description of a Toda-type tt*-structure and we see that the anti-symmetry condition for \(l=0, 1\) are essential. From Theorem \ref{prop3.1} and Lemma \ref{lem3.1}, we classify Toda-type tt*-structures as follows.
\begin{cor}\label{cor3.1}
	Let \(\sim_w\) be an equivalence relation on solutions to the tt*-Toda equation defined by
	\begin{equation}
		(w^2_0,\cdots,w^2_n) \sim_w (w_0^1,\cdots,w_n^1)\ \ \ \Longleftrightarrow\ \ \ \exists s \in \mathbb{Z}\ \ \ {\rm s.t.}\ \ \ w_j^2= w_{j+s}^1. \nonumber
	\end{equation}
	Then, there is an one-to-one correspondence between
	\begin{itemize}
		\item [(1)] an isomorphism class of a Toda-type tt*-structure,
		
		\item [(2)] a \(\mathbb{Z}_{n+1}\)-fixed point with a non-vanishing \(\Phi\),
		
		\item [(3)] an equivalence class of a solution to the tt*-Toda equation.
	\end{itemize}
\end{cor}
\begin{proof}
	\((1) \rightarrow (2):\) It follows from Lemma \ref{lem3.1}.\vspace{2mm}
	
	\noindent \((2) \rightarrow (3):\) Let \((E_1,\eta_1,g_1,\Phi_1), (E_2,\eta_2,g_2,\Phi_2)\) be representatives of a \(\mathbb{Z}_{n+1}\)-fixed point. Then there exists an isomorphism \(\mathcal{T}:(E_2,\eta_2,g_2,\Phi_2) \stackrel{\sim}{\longrightarrow} (E_1,\eta_1,g_1,\Phi_1)\) of tt*-structures. Choose the frame \(\tau^a= (\tau_0^a,\cdots,\tau_n^a)\) of \(E_a\) such that \(\Phi(\tau_j^a) = \omega^j \tau_j^a dt^a\), where \(t^a\) is the coordinate of \(\mathbb{C}^*\), and \(\eta_a(\tau_i^a,\tau_j^a) = \delta_{ij}\ (a=1,2)\). Since \(\{\omega^j dt^1\}_{j=0}^{n}\) and \(\{\omega^j dt^2\}_{j=0}^{n}\) are eigenvalues of \(\Phi\), we have \(\frac{dt^2}{dt^1} = \omega^k\ (k \in \mathbb{Z})\). Put \(\mathcal{T}(\tau^2) = \tau^1 \cdot T\). From \(\eta_1(\tau_i^1,\tau_j^1) = \eta_2(\tau_i^2,\tau_j^2) = \delta_{ij}\), we have
	\begin{equation}
		T = \left(\begin{array}{ccc}
			a_0 & & \\
			& \ddots & \\
			& & a_n
		\end{array}\right)M,\ \ \ a_j^2=1, \nonumber
	\end{equation}
	where \(M\) is a permutation matrix such that
	\begin{equation}
		\omega^k \left(\begin{array}{ccc}
			1 & & \\
			& \ddots & \\
			& & \omega^n
		\end{array}\right) = M^{-1}\left(\begin{array}{ccc}
			1 & & \\
			& \ddots & \\
			& & \omega^n
		\end{array}\right)M. \nonumber
	\end{equation}
	By replacing \(\tau^2\) by \(\tau^2 \cdot M^{-1}\), without loss of generality we can assume that \(t^2=t^1\) and \(k=0\). Let \((w_0^j,\cdots,w_n^j)\) be solutions to the tt*-Toda equation given by \((E_a,\eta_a,g_a,\Phi_a)\ (a=1,2)\)
	and \(L\) be the matrix in the proof of Theorem \ref{prop3.1}, then we have
	\begin{equation}
		\left(\begin{array}{ccc}
			e^{w^2_0} & & \\
			& \ddots & \\
			& & e^{w^2_n}
		\end{array}\right) = (L^{-1}TL)^{-1}\left(\begin{array}{ccc}
			e^{w^1_0} & & \\
			& \ddots & \\
			& & e^{w^1_n}
		\end{array}\right)(L^{-1}TL), \nonumber
	\end{equation}
	and then, there exists a \(s \in \mathbb{Z}\) such that \(w_j^2 = w_{j+s}^1\) for \(j=0,\cdots,n\).\vspace{2mm}
	
	\noindent \((3) \rightarrow (1)\): Let \(\{w_j^1\}_{j=0}^{{n}}\) and \(\{w_j^2\}_{j=0}^{n}\) be solutions to the tt*-Toda equation such that \(w_j^2=w_{j+s}^1\ (s \in \mathbb{Z})\) and \((E_a,\eta_a,g_a,\Phi_a)\ (a=1,2)\) tt*-structures with the frame \(\{e_j^a\}_{j=0}^{n}\) such that
	\begin{equation}
		\eta(e_i^a,e_j^a) = \left\{\begin{array}{ll}
			\delta_{i,l-1-j} & {\rm if}\ 0 \le a \le l-1,\\
			\delta_{i,n+l-j} & {\rm if}\ l \le a \le n,
		\end{array}\right. ,\ \ \ \ \ g(e_i^a,e_j^a) = e^{w_j^a}\delta_{ij}, \nonumber
	\end{equation}
	\begin{equation}
		\Phi (e_0^a,\cdots,e_n^a) = (e_0^a,\cdots,e_n^a) \left(\begin{array}{cccc}
			& & & 1 \\
			1 & & &   \\
			& \ddots & &  \\
			& & 1 & 
		\end{array} \right)dt,\ \ \ t \in \mathbb{C}^*. \nonumber
	\end{equation} 
	We define a map \(\mathcal{T}:E_2 \rightarrow E_1\) by
	\begin{equation}
		\mathcal{T}(e_0^2,\cdots,e_n^2) = (e_0^1,\cdots,e_n^1)\left(\begin{array}{cccc}
			& & & 1 \\
			1 & & & \\
			& \ddots & & \\
			& & 1 &
		\end{array}\right)^s, \nonumber
	\end{equation}
	then \(\mathcal{T}:(E_2,\eta_2,g_2,\Phi_2) \rightarrow (E_1,\eta_1,g_1,\Phi_1)\) is an isomorphism of tt*-structures.
\end{proof}

From Corollary \ref{cor3.1}, we can see that any solutions to the tt*-Toda equation are equivalent to a solution to the tt*-Toda equation with the anti-symmetry condition for \(l=0\) or \(l=1\). Hence, the tt*-equation is characterized by the \(l=0,1\) in the anti-symmetry condition. In the following section, we give some observations of the intrinsic description of the tt*-Toda equation and the classification.

\section{Applications}\label{s4}
In section \ref{s3}, we obtain an intrinsic description of the tt*-Toda equation and we classify Toda-type tt*-structures
\begin{align}
	&\{\text{isomorphic class of a Toda-type tt*-structure}\} \nonumber\\
	&\ \ \ \ \ \ \ \ \ \ \simeq \{\mathbb{Z}_{n+1}\text{-fixed point with a non-vanishing \(\Phi\)}\}. \nonumber
\end{align}
In this section, we give some observations from the intrinsic description and the classification. In section \ref{s4}, we explain the equivalence of the anti-symmetry condition occurred from the classification of Toda-type tt*-structures. As an example, we consider the tt*-Toda equation which involves two unknown functions introduced by Guest and Lin \cite{GL2014}. In section \ref{s4-2}, we review the results of Fredrickson and Neitzke \cite{FN2021} and we give the relation between radial solutions to the tt*-Toda equation and representations of the \(W\)-algebra of \(\mathfrak{sl}_{n+1} \mathbb{C}\).

\subsection{An equivalence of \(l\) in the anti-symmetry condition}\label{s4-1}
We described the tt*-Toda equation with the anti-symmetry condition as a \(\mathbb{Z}_{n+1}\)-fixed point. In Theorem \ref{prop3.1}, the \(l\) in the anti-symmetry condition
\begin{equation}
	\left\{\begin{array}{l}
		w_0 + w_{l-1} = 0, w_1 + w_{l-2} = 0, \cdots \\
		w_l + w_n=0, w_{l+1} + w_{n-1} = 0, \cdots.
	\end{array}\right.\ \ \ \ \ l \in \mathbb{Z}, \nonumber
\end{equation}
can take on any of the two values \(0,1\). From Corollary \ref{cor3.1}, the tt*-Toda equation with the anti-symmetry condition for \(l \in \mathbb{Z}\) is equivalent to the tt*-Toda equation with the condition
\begin{equation}
	w_j + w_{n-j} = 0,\ \  \ \ j = 0,\cdots,n, \nonumber
\end{equation}
or
\begin{equation}
	w_0=0,\ \ \ w_j + w_{n+1-j} = 0,\ \  \ \ j = 1,\cdots,n. \nonumber
\end{equation}
It can also be proved more directly.
\begin{thm}\label{thm1.2}
	Let \(\{w_j\}_{j=0}^{n}\) be solutions to the tt*-Toda equation satisfying the anti-symmetry condition with \(l\)
	\begin{equation}
		\left\{\begin{array}{l}
			w_0 + w_{l-1} = 0, w_1 + w_{l-2} = 0, \cdots \\
			w_l + w_n = 0, w_{l+1} + w_{n-1} = 0, \cdots.
		\end{array}\right. \nonumber
	\end{equation}
	Then, there exists a solution \(\{\tilde{w}_j\}_{j=0}^{n}\) to the tt*-Toda equation satisfying the anti-symmetry condition
	\begin{equation}
		\left\{\begin{array}{cc}
			\tilde{w}_j + \tilde{w}_{n-j} = 0 & \text{if l or n is an even number,} \\
			\tilde{w}_j + \tilde{w}_{n+1-j} = 0 & \text{otherwise}
		\end{array}
		\right. \nonumber
	\end{equation}
	for \(j=0,\cdots,n\), and
	\begin{equation}
		(w_0,\cdots,w_n) \sim_w (\tilde{w}_0,\cdots,\tilde{w}_n), \nonumber
	\end{equation}
	where \(\sim_w\) is the equivalence relation defined in Corollary \ref{cor3.1}.
\end{thm}
\begin{proof}
	When \(l=2m\ (m \in \mathbb{Z})\), we put \(\{\tilde{w}_j\} = \{w_{j+l/2}\}\) then we have
	\begin{equation}
		\tilde{w}_j = w_{j+m} = -w_{m-1-j} = -\tilde{w}_{-1-j} = -\tilde{w}_{n-j}. \nonumber
	\end{equation}
	When \(l=2m+1\ (m \in \mathbb{Z})\), we put \(\{\tilde{w}_j\} = \{w_{j+(l-1)/2}\}\) then we have
	\begin{equation}
		\tilde{w}_j = w_{j+m} = -w_{m-j} = -\tilde{w}_{-j} = -\tilde{w}_{n+1-j}. \nonumber
	\end{equation}
	When \(l=2m+1, n=2s (m,s \in \mathbb{Z})\), we put \(\{\tilde{w}_j\} = \{w_{j-(n-l+1)/2}\}\) then we have 
	\begin{equation}
		\tilde{w}_j = w_{j-(s-m)} = -w_{s+m-j} = -\tilde{w}_{n-j}. \nonumber
	\end{equation}
\end{proof}
Thus, Theorem \ref{thm1.2} gives the equivalence of the cases \(l=0,\cdots,n\) in the anti-symmetry condition. Every \(\mathbb{Z}_{n+1}\)-fixed point is represented by a tt*-structure whose tt*-equation is 
the tt*-Toda equation with the anti-symmetry condition \(w_j+w_{n-j}=0\) or \(w_j+w_{n+1-j}=0\). From the intrinsic point of view, we can identify the case \(l \in \mathbb{Z}\) with the case \(l=0\) or \(l=1\). In particular, if \(n\) is an even number, the condition \(w_j+w_{n-j} = 0\) is essential.\vskip\baselineskip

We apply Theorem \ref{thm1.2} to the tt*-Toda equation which involves two unknown functions.
\begin{ex}
	In \cite{GL2014}, Guest and Lin proved (Corollary 2.3 of \cite{GL2014}) that the tt*-Toda equation which involves two unknown functions can be written of the form
	\begin{equation}
		\left\{
		\begin{array}{l}
			w_{z\overline{z}} = e^{aw} - e^{v-w},\\
			v_{z\overline{z}} = e^{v-w} - e^{-bv},
		\end{array}
		\right. \nonumber
	\end{equation}
	where \(a,b \in \{1,2\}\), and there are ten possibilities for \(l\) in the anti-symmetry condition (see Table \ref{table1})
	\begin{equation}
		\left\{\begin{array}{l}
			w_0 + w_{l-1} = 0, w_1 + w_{l-2} = 0, \cdots \\
			w_l + w_n = 0, w_{l+1} + w_{n-1} = 0, \cdots.
		\end{array}\right. \nonumber
	\end{equation}
	\begin{table}[h]
		\centering
		\begingroup
		\renewcommand{\arraystretch}{1.5}
		\caption{Table 3 of \cite{GL2014}}\label{table1}
			\begin{tabular}{c|cccccccccc}
				\(n\) & \(3\) & \(4\) & \(4\) & \(5\) & \(3\) & \(4\) & \(4\) & \(5\) & \(4\) & \(5\) \\
				\(l\) & \(0\) & \(0\) & \(1\) & \(1\) & \(2\) & \(2\) & \(3\) & \(3\) & \(4\) & \(5\) \\ \hline
				\(w\) & \(w_0\) & \(w_0\) & \(w_1\) & \(w_1\) & \(w_3\) & \(w_4\) & \(w_4\) & \(w_5\) & \(w_0\) & \(w_0\) \\
				\(v\) & \(w_1\) & \(w_1\) & \(w_2\) & \(w_2\) & \(w_0\) & \(w_0\) & \(w_0\) & \(w_0\) & \(w_1\) & \(w_1\) \\ \hline
				\(a\) & \(2\) & \(1\) & \(2\) & \(1\) & \(2\) & \(2\) & \(1\) & \(1\) & \(2\) & \(1\) \\
				\(b\) & \(2\) & \(2\) & \(1\) & \(1\) & \(2\) & \(1\) & \(2\) & \(1\) & \(1\) & \(1\) \\
			\end{tabular}
		\endgroup \nonumber
	\end{table}
	
	From Theorem \ref{thm1.2}, the possibilities for \(n,l\) can be reduced to just three cases \((n,l) = (3,0), (4,0), (5,1)\). Two pairs \((n,l)\)'s are equivalent if and only if the corresponding \((a,b)\)'s are equal up to permutations.\vskip\baselineskip
	
	Guest and Lin reduced ten cases to four cases \((a,b) = (2,2), (1,2), (2,1), (1,1)\). On the other hand, we reduce three cases. This is because we obtain the case \((a,b) = (2,1)\) from the case \((a,b) = (1,2)\) by transforming \((w,v)\) into \((-v,-w)\). \qed
\end{ex}

From Theorem \ref{thm1.2}, we obtain a generalization of Corollary 2.3 of \cite{GL2014}.
\begin{cor}\label{cor4.2}
	Any system arising from the tt*-Toda equation with the anti-symmetry condition for \(l\) which involves \(m\) unknown functions can be written in the form (see Table 2) 
	\begin{equation}
		\left\{
		\begin{array}{l}
			(v_0)_{t \overline{t}} = e^{av_0} - e^{v_1-v_0}, \\
			(v_j)_{t \overline{t}} = e^{v_j - v_{j-1}} - e^{v_{j+1} - v_j},\ \ \ \ \ \ \ \ \ \ 1 \le j \le m-2, \\
			(v_{m-1})_{t \overline{t}} = e^{v_{m-1} - v_{m-2}} - e^{bv_{m-1}},
		\end{array}
		\right. \nonumber
	\end{equation}
	where \(a,b \in \{1,2\}\).
\end{cor}
\begin{proof}
	\begin{table}[h]
		\centering
		\begingroup
		\renewcommand{\arraystretch}{1.5}
		\caption{The case \(l=0,1\)}\label{table2}
			\begin{tabular}{c|cccc}
				\(n\) & \(2m-1\) & \(2m\) & \(2m\) & \(2m+1\) \\
				\(l\) & \(0\) & \(0\) & \(1\) & \(1\) \\ \hline
				\(v_0\) & \(w_0\) & \(w_0\) & \(w_0\) & \(w_1\) \\
				\(\vdots\) & \(\vdots\) & \(\vdots\) & \(\vdots\) & \(\vdots\) \\
				\(v_{m-1}\) & \(w_{m-1}\) & \(w_{m-1}\) & \(w_{m-1}\) & \(w_m\) \\ \hline
				\(a\) & \(2\) & \(2\) & \(2\) & \(1\) \\
				\(b\) & \(2\) & \(1\) & \(1\) & \(1\) \\
			\end{tabular}
		\endgroup \nonumber
	\end{table}
	From Theorem \ref{thm1.2}, there are \(4m+2\) possibilities for \((n,l)\), but it is enough to consider the case \(l=0,1\). We summarize in Table \ref{table2}.
\end{proof}

Hence, the tt*-Toda equation which involves \(m\) unknown functions can be written by any of just three distinct p.d.e. systems.

\subsection{\(\rho\)-actions on the moduli of wild Higgs bundles}\label{s4-2}
In Theorem 5.3 of \cite{FN2021}, Fredrickson and Neitzke suggested that the effective Virasoro central charge of a representations of \(W\)-algebra of \(\mathfrak{sl}{_{n+1} \mathbb{C}}\) can be described by a \(L^2\)-norm \(\mu\) of a fixed point of the  \(\rho\)-action (Definition \ref{def2.5}) in the moduli of wild Higgs bundles (see definition in section 2 of \cite{FN2021}) over \(\mathbb{C}P^1\). The \(L^2\)-norm is defined by the harmonic metric of the fixed point (see definition in (4.2) of \cite{FN2021}) and it is an extension of a real-valued function on the moduli of harmonic bundles over a compact Riemann surface introduced by Hitchin (\cite{H1987}). Let \(p,q \in \mathbb{Z}\) be coprime, then irreducible representations of \(W\)-algebra \(\mathcal{W}_{n+1}\) parameterized by dominant weights of \(\mathfrak{sl}_{n+1} \mathbb{C}\) (\cite{BS1993} and Theorem 5.1 of \cite{GO2022}) describe a \(2\)-dimensional chiral conformal field theory called \((p,q)\)-minimal model. A direct sum of finitely many irreducible representations of \(\mathcal{W}_{n+1}\) is the Hilbert space of \((p,q)\)-minimal model. In particular, Fredrickson and Neitzke observed the case \((p,q) = (n+1,N)\).\vskip\baselineskip

For a fixed point of the  \(\rho\)-action in the moduli of wild Higgs bundles over \(\mathbb{C}P^1\), the Hitchin equation can be written as
\begin{equation}
	\frac{1}{4}\left(\frac{d^2w_j}{dr^2} + \frac{1}{r}\frac{dw_j}{dr}\right) = e^{w_j-w_{j-1}} - e^{w_{j+1}-w_j}, \nonumber
\end{equation}
in the radial coordinate \(r = |t|\). Since the \(\{w_j\}\) does not satisfy the anti-symmetry condition, this differential equation is not the tt*-Toda equation. Regarding the radial tt*-Toda equation, Guest, Its and Lin found all global radial solutions to the tt*-Toda equation \cite{GIL20151}, \cite{GIL20152} and Guest and Otofuji gave a relation between radial solutions to the tt*-Toda equation with the anti-symmetry condition \(l=0\) and representation theory \cite{GO2022}).\vskip\baselineskip

In this section, we observe a relation between radial solutions to the tt*-Toda equation with the anti-symmetry condition \(l \in \mathbb{Z}\) and representations of the \(W\)-algebra of \(\mathfrak{sl}_{n+1} \mathbb{C}\). First, we consider the \(\rho\)-action (Definition \ref{def2.5}) on Toda-type tt*-structures and we see that fixed points of the \(\rho\)-action in the set of Toda-like tt*-structures give the radial tt*-Toda equations. We follow the argument of section 3.4 in \cite{FN2021}.

\begin{prop}\label{prop4.2}
	Let \((E,\eta,g,\Phi)\) be a tt*-structure of rank \(n+1\) over \(\mathbb{C}^*\) with a non-vanishing \(\Phi\). If \((E,\eta,g,\Phi)\) is a representative of \(\mathbb{Z}_{n+1}\)-fixed point and there exists an isomorphism
	\begin{equation}
		\mathcal{T}_z: (\rho_z^*E,\rho_z^*\eta,\rho_z^*g,e^{-z}\rho_z^*\Phi) \stackrel{\sim}{\longrightarrow} (E,\eta,g,\Phi), \nonumber
	\end{equation}
	of tt*-structures for all \(z \in \sqrt{-1}\mathbb{R}\) and \(T_z(a)\) is continuous in \(z\) for all \(a \in E\), then the tt*-equation of \((E,\eta,g,\Phi)\) is equivalent to the tt*-Toda equation with the radial condition.\\
	We call the isomorphism class of the tt*-structure \((E,\eta,g,\Phi)\) above the \(\rho\)-fixed point.
\end{prop}
\begin{proof}
	Since \((E,\eta,g,\Phi)\) is a representative of a \(\mathbb{Z}_{n+1}\)-fixed point, there exists a local frame \(\tau = (\tau_0,\cdots,\tau_n)\) of \(E\) such that \(\eta(\tau_i,\tau_j) = \delta_{ij}\) and \(\Phi(\tau_j) = \omega^j \tau_j dt\). Put \(\mathcal{T}_z(t,\tau_j(\rho_z(t))) = \sum_{i=0}^{n}T^z_{ij}(t)\tau_i(t)\), then from the assumption we have \(T_{ij}^z\) is continuous in \(z\). Since
	\begin{equation}
		\Phi_{\frac{d}{dt}}(\mathcal{T}_z(t,\tau_j(\rho_z(t)))) = \mathcal{T}_z((e^{-z}\rho_z^* \Phi)_{\frac{d}{dt}}(t,\tau_j(\rho_z(t)))), \nonumber
	\end{equation}
	and
	\begin{equation}
		(\rho_z^* \eta)_t(t,\tau_i(\rho_z(t)),(t,\tau_j(\rho_z(t)))) = \eta_t(\mathcal{T}_z(t,\tau_i(\rho_z(t))),\mathcal{T}^z(t,\tau_J(\rho_z(t)))), \nonumber
	\end{equation}
	we have \(T^z_{ij} = \varepsilon_j \delta_{ij}\), where \(\varepsilon_j \in \{-1,1\}\). Since \(T_{ij}^0 = \delta_{ij}\) and \(T^z_{ij}\) is continuous in \(z\), we have \(\varepsilon_j = 1\) and then \(\mathcal{T}^z = Id\). Thus, we obtain 
	\begin{align}
		g(\tau_i,\tau_j) \circ \rho_z(t) 
		&= (\rho_z^* g)_t((t,\tau_i(\rho_z(t))),(t,\tau_j(\rho_z(t)))) \nonumber\\
		&= g_t(\mathcal{T}^z(t,\tau_i(\rho_z(t))),\mathcal{T}^z(t,\tau_j(\rho_z(t)))) \nonumber\\
		&= g(\tau_i,\tau_j)(t). \nonumber
	\end{align}
	From the proof of Theorem \ref{prop3.1}, we obtain the stated result.
\end{proof}

From Proposition \ref{prop4.2} and Corollary \ref{cor3.1}, we obtain the following relation.
\begin{cor}
	There is an one-to-one correspondence between
	\begin{itemize}
		\item [(I)] a \(\mathbb{Z}_{n+1}\)-fixed point and \(\rho\)-fixed point with a non-vanishing \(\Phi\),
		
		\item [(II)] an equivalence class of a solution to the tt*-Toda equation with the radial condition.
	\end{itemize}
	We call a representative of a \(\mathbb{Z}_{n+1}\)-fixed point and \(\rho\)-fixed point with a non-vanishing \(\Phi\) a radial Toda-type tt*-structure.
\end{cor}\vskip\baselineskip

Next, we characterize radial Toda-type tt*-structures by using the result introduced by Guest, Its and Lin. In \cite{GIL20151}, \cite{GIL20152}, \cite{GIL2020}, Guest, Its and Lin showed that radial solution \(\{w_j\}_{j=0}^{n}\) to the tt*-Toda equation is characterized by the asymptotic behaviour
\begin{equation}
	w_j \sim -m_j \log{|t|}\ \ \ \ \ {\rm as}\ \ \ t \rightarrow 0, \nonumber
\end{equation}
where \(m_j \in \mathbb{R}\) and \(m_{j-1} - m_j +2 \ge 0\). By using this characterization and the \(l\) in the anti-symmetry condition, we characterize Toda-type tt*-structures.

\begin{prop}\label{prop4.1}
	Let \(\sim_l\) be an equivalence relation on \((m_0,\cdots,m_n) \in \mathbb{R}\) defined by
	\begin{equation}
		(\tilde{m}_0,\cdots,\tilde{m}_n) \sim_l (m_0,\cdots,m_n)\ \ \ \Longleftrightarrow\ \ \ \exists k \in \mathbb{Z}\ \ {\rm s.t.}\ \ \tilde{m}_j = m_{j+k}, \nonumber
	\end{equation}
	where \(m_{n+1+a} = m_a\), and we denote the equivalence class by \([(m_0,\cdots,m_n)]\). Put
	\begin{equation}
		\mathfrak{M}_l = \left\{(l,[m_0,\cdots,m_n])\ | \ m_j + m_{n+l-j} = 0,\ m_{j-1} - m_j + 2 \ge 0\right\},\ \ \ l = 0,1, \nonumber
	\end{equation}
	where \(m_{n+1}=m_0,m_{-1}=m_n\), then if \(n\) is an even number there is a bijection between
	\begin{equation}
		\mathfrak{M}_0 \simeq \left\{\text{isomorphic class of a radial Toda-type tt*-structure}\right\}, \nonumber
	\end{equation}
	\noindent if \(n\) is an odd number there is a bijection between
	\begin{equation}
		\mathfrak{M}_0 \cup \mathfrak{M}_1 \simeq \left\{\text{isomorphic class of a radial Toda-type tt*-structure}\right\}, \nonumber
	\end{equation}\vspace{2mm}
	In this paper, we call \((l,[m_0,\cdots,m_n])\) the \(l\)-asymptotic data.
\end{prop}
\begin{proof}
	It follows from Theorem \ref{thm3.2}, Proposition \ref{prop4.2} and Theorem 2.1 of \cite{GO2022}
\end{proof}
Hence, a radial Toda-type tt*-structure is characterized by the \(l\)-asymptotic data \((l,\left[m_0,\cdots,m_n\right])\). We use the \(l\)-asymptotic data and we apply the observation of Fredrickson, Neitzke to radial Toda-type tt*-structures.\vskip\baselineskip

In Theorem 5.3 of \cite{FN2021}, Fredrickson and Neitzke observed a relation between \(\rho\)-fixed points of the moduli of wild Higgs bundles over \(\mathbb{C}P^1\) and the effective Virasoro central charge of the representations of \(W\)-algebra of \(\mathfrak{sl}{_{n+1} \mathbb{C}}\). The \(L^2\)-norm \(\mu\) of a \(\rho\)-fixed point in the harmonic metric is related to the effective Virasoro central charge by a simple formula (Theorem 5.3 of \cite{FN2021}). In section 5 of \cite{GO2022}, Guest and Otofuji gave a relation between solutions to the tt*-Toda equation with the anti-symmetry condition \(l = 0\) and the effective Virasoro central charge of the representations of \(W\)-algebra of \(\mathfrak{sl}{_{n+1} \mathbb{C}}\). The asymptotic behaviour of solutions is related to the effective Virasoro central charge by formula (5.5) of \cite{GO2022}. In both cases, the effective Virasoro central charge can be described by the asymptotics of solutions to the p.d.e. (the Hitchin equation in \cite{FN2021}, the tt*-Toda equation in \cite{GO2022}).\vskip\baselineskip

In this paper, we give a similar observation for the case of radial Toda-type tt*-structures and we associate solutions to the tt*-Toda equation (satisfying the anti-symmetry condition \(l \in \mathbb{Z}\)) with the effective Virasoro central charge of \((n+1,N)\)-minimal model of the \(W\)-algebra \(\mathcal{W}_{n+1}\). Here, the asymptotic data \(m_j\)'s satisfies the anti-symmetry condition for \(l\) and the corresponding effective Virasoro central charge \(c_{\rm eff}\) is determined uniquely by an equivalence class \((l,[m_0,\cdots,m_n])\).

\begin{cor}\label{cor6.6}
	Let \((E,\eta,g,\Phi)\) be a representative of a radial Toda-type tt*-structure with the \(l\)-asymptotic data \((l,[m_0,\cdots,m_n])\).\\
	If \((m_0,\cdots,m_n) \in \mathbb{Q}^{n+1}\) and \(m_{j-1} + m_j + 2 > 0\) for all \(j\), then there exist \(N \in \mathbb{N}\) and the highest weight representation \(\Lambda\) in the \((n+1,N)\)-minimal model of \(\mathcal{W}_{n+1}\) such that \(\Lambda\) is parameterized by
	\begin{equation}
		[m_n-m_0,m_0-m_1,\cdots,m_{n-1}-m_n], \nonumber
	\end{equation}
	and the effective Virasoro central charge \(c_{\rm eff}\) with \(\Lambda\) is given by
	\begin{equation}
		c_{\rm eff} = n - \frac{3(N + n + 1)}{n+1}\sum_{j=0}^{n}m_j^2, \nonumber
	\end{equation} 
	where \((l,[m_0,\cdots,m_n])\) is the \(l\)-asymptotic data of \((E,\eta,g,\Phi)\).
\end{cor}
\begin{proof}
	For the definition of the \(W\)-algebra \(\mathcal{W}_{n+1}\) of \(\mathfrak{sl}_{n+1} \mathbb{C}\), we refer to section 5 of \cite{FN2021} and \cite{BS1993}. From Proposition 5.1 of \cite{FN2021}, the highest representation of \(\mathcal{W}_{n+1}\) is determined by a \(n+1\)-tuple \((b_0,\cdots,b_n) \in \mathbb{Z}_{\ge 0}^{n+1}\). Let \(m_j = 2P_j/Q\ (P_j \in \mathbb{Z}, Q \in \mathbb{Z}_{\ge 2})\), then we put \(N = (n+1)(Q-1)\) and
	\begin{equation}
		b_j = \frac{N + n + 1}{2(n+1)}(m_{j-1}-m_j+2) - 1 \in \mathbb{Z}_{\ge 0},\ \ \ \ \ j=0,\cdots,n. \nonumber
	\end{equation}
	From Proposition 4.1 and Theorem 5.3 of \cite{FN2021}, we have
	\begin{equation}
		\frac{N+n+1}{4(n+1)}\sum_{j=0}^{n}m_j^2 = \frac{1}{12}(n - c_{\rm eff}). \nonumber
	\end{equation}
	Hence, we obtain
	\begin{equation}
		c_{\rm eff} = n - \frac{3(N + n + 1)}{n+1}\sum_{j=0}^{n}m_j^2. \nonumber
	\end{equation}
\end{proof}

\begin{rem}
	When we choose \(l=0\), the asymptotic data \((0,[m_0,\cdots,m_n])\) satisfies
	\begin{equation}
		m_j + m_{n-j} = 0,\ \ \ \ \ m_{j-1} - m_j + 2 \ge 0. \nonumber
	\end{equation}
	In this case (i.e. \(l = 0\) case), the \(l\)-asymptotic data \((0,[m_0,\cdots,m_n])\) is equal to the asymptotic data introduced by Guest, Its, Lin (\cite{GIL20151}, \cite{GIL20152}, \cite{GIL2020}). Then, the relation between solutions to the radial tt*-Toda equation and the \(W\)-algebra was given by Guest, Otofuji (section 5 of \cite{GO2022}). \qed
\end{rem}

\section*{Acknowledgement}
The author would like to thank Professor Martin Guest for his considerable support and thank Professor Takashi Otofuji for useful conversations. This paper is a part of the outcome of research performed under a Waseda University Grant for Special Research Projects (Project number: 2025C-104).

\section*{Conflict of interests}
The author has no conflicts to disclose.

\bibliography{mybibfile}
\bibliographystyle{plain}
	
	\em
	\noindent
	Department of Applied Mathematics\newline
	Faculty of Science and Engineering\newline
	Waseda University\newline
	3-4-1 Okubo, Shinjuku, Tokyo 169-8555\newline
	JAPAN

\end{document}